# Abel-Jacobi theorem

Seddik Gmira

USMBA

## 1 Introduction

The Abel-Jacobi theorem is an important result of algebraic geometry. The theory of divisors and the Riemann bilinear relations are fundamental to the developement of this result: if a point $O$ is fixed in a Riemann compact surface $X$ of genus $g$, the Abel-Jaobi map identifies the Picard group $Pic_O(X)$ the quotient of divisors of a group of degree zero on the sub-group of divisors associated to meromorphic functions. The Riemann surface of genus $g \geq 1$ can be embedded in the Jacobian variety $Jac(X)$ via the Abel-Jacobi. In fact we generally have a map:

$$X^{(g)} = X^g / \mathfrak{S}_g \longrightarrow Jac(X)$$

such that $X^{(g)}$ may be provided with an analytical structure. Indeed the two sets $X^{(g)} = X^g / \mathfrak{S}_g$, $Jac(X)$ are algebraic varieties and the map

$$X^{(g)} \longrightarrow Jac(X)$$

is surjective. For reasons of dimension we can verify that is finite fibers. In fact this is a birational map.

## 2 Riemann bilinear relations

Let $X$ be a compact Riemannian surface. Recalling that,

$$H_1(X, \mathbb{Z}) \approxeq \mathbb{Z}^{2g} \text{ and } H^1_{dR}(X, \mathbb{R}) \approxeq \mathbb{R}^{2g}$$

where $g$ is the genus of $S$. The following map

$$\begin{array}{ccc} H_1(X, \mathbb{Z}) \times H^1_{dR}(X, \mathbb{R}) & \longrightarrow & \mathbb{R} \\ (\gamma, \omega) & \longrightarrow & \int_\gamma \omega \end{array}$$

makes these two spaces in duality: for a basis $(\gamma_1, ..., \gamma_{2g})$ in $H_1(X, \mathbb{Z})$ there exist a dual basis
$(\omega_1, ..., \omega_{2g})$ in $H^1_{dR}(X, \mathbb{R})$ such that for $i, j = 1, ..., 2g$

$$\int_{\gamma_i} \omega_j = \delta_{ij}$$

The intersection product

$$H_1(X,\mathbb{Z}) \times H_1(X,\mathbb{Z}) \longrightarrow \mathbb{Z}$$
$$(\gamma_1, \gamma_2) \longrightarrow \gamma_1 \# \gamma_2$$

defines an antisymetric bilinear form on $H_1(X,\mathbb{Z})$, which has a corresponding symplectic bases

**Proposition 1** *For any symplectic basis $(a_1, ..., a_g, b_1, ..., b_g)$ of $H_1(X,\mathbb{Z})$ and for any closed 1-formes $\eta$ and $\eta'$ on the surface $X$ we have*

$$\int_X \eta \wedge \eta' = \sum_{k=1}^{g} \left( \int_{a_i} \eta \int_{b_i} \eta' - \int_{a_i} \eta' \int_{b_i} \eta \right)$$

**Preuve.** Let $(a_1, ..., a_g, b_1, ..., b_g)$ be a symplectic basis of $H_1(X,\mathbb{Z})$ associated with a cuuting $S$ into a $4g$−Gones quotes $\Delta$: $A_1 B_1 A'_1 B'_1, ....., A_g B_g A'_g B'_g$, where $A_i$ and $A'_i$ are identified by the map $\varphi_i$ and $B_i$, $B'_i$ are identified by the map $\psi_i$ as in the following figure. Differential formes can be seen as differential formes on $\Delta$. Since this last is simply connected, so there exist a function $f$ such that $df = \eta$. So for each $x \in A$ and for each $y \in B$ we have:

$$(1) : \int_{b_i(x)} df = \int_{b_i} \eta = f \circ \varphi_i(x) - f(x)$$

$$(2) : \int_{a_i(x)} df = \int_{a_i} \eta = f(x) - f \circ \psi_i(x)$$

Stokes formula implies

$$\begin{aligned}
\int_S \eta \wedge \eta' &= \int_\Delta \eta \wedge \eta' \\
&= \int_D d(f\eta') \\
&= \int_\Delta f\eta' \\
&= \sum_{k=1}^g \int_{A_i+B_i-A'_i-B'_i} f\eta'
\end{aligned}$$

And it follows from the formulas (1) and (2):

$$\int_{A_i-A'_i} f\eta' = \int_{A_i} (f - f \circ \varphi_i(x))\,\eta' = -\int_{b_i} \eta \int_{a_i} \eta'$$

$$\int_{B_i-B'_i} f\eta' = \int_{B_i} (f - f \circ \psi_i(x))\,\eta' = \int_{a_i} \eta \int_{b_i} \eta'$$

which proves equality ∎

**Remarque 2** *If the surface $X$ is provided with a riemann structure, and if $\eta\ \eta'$ are holomorphic 1-forms, then $\int_X \eta \wedge \eta' = 0$*

**Proposition 3** *Let $X$ be a compact Riemannian of which is fixed $2g$ simple closed curves $(a_1, ..., a_g, b_1, ..., b_g)$, forming a symplectic basis of the space $H_1(X, \mathbb{Z})$ and let $\omega_1$ be a holomorphic 1-form on $X$ and $\omega_2$ non-sigular 1-meromorphic form along all the curves $a_i\ b_i$. Given a point $z_0 \in X - \{a_i\ b_i\}$ such that, $u(z) = \int_{z_0}^z \omega_1$, then*

$$2i\pi \sum \mathrm{Res}(u.\omega_2) = \sum_{i=1}^g \left( \int_{a_i} \omega_1 \int_{b_i} \omega_2 - \int_{a_i} \omega_2 \int_{b_i} \omega_1 \right)$$

**Preuve.** The proposal follows from the Residue formula and equations (1) and (2) : $2i\pi \sum \mathrm{Res}(u.\omega_2) = \int_{\partial \Delta} u.\omega_2$ ∎

Whether now $(a_1, ..., a_g, b_1, ..., b_g)$ is a $2g$ simple closed curves on a compact Riemann surface $X$ which form basis of the space $H_1(X, \mathbb{Z})$ and $(\omega_1, ..., \omega_g)$ is a fixed basis of the space of 1-holomorphic forms on $X$.

**Dfinition 4** *Let's call the period matrices $A, B \in \mathfrak{M}_g(\mathbb{C})$ defined by*

$$A_{ij} = \int_{a_i} \omega_j$$

$$B_{ij} = \int_{b_i} \omega_j$$

**Thorme 5** *(Riemann bilinear relations)*

1. The matrix $A$ is invertible

2. The matrix $\Omega = A^{-1}B$ is symetrical and its imaginary part

$$\operatorname{Im}\Omega = (\operatorname{Im}\Omega_{ij})_{i,j \leq g}$$

   is positive definite

**Preuve.** Whether $\lambda = (\lambda_1, .., \lambda_g) \in \mathbb{C}^g$ such that $\sum_{i=1}^{g} \lambda_i A_{ij} : j = 1, ..., g$. Consider the holomorphic 1-form

$$\omega = \sum_{i=1}^{g} \lambda_i \omega_i$$

By definition of the matrix $A$, we have:

$$\int_{a_i} \omega = 0 = \sum_{i=1}^{g} \lambda_i A_{ij}$$

so is

$$\int_{a_i} \overline{\omega} = 0$$

Then it follows from the Proposition1,

$$\int_{a_i} \omega \wedge \overline{\omega} = 0 : \omega = 0$$

so $\lambda_i = 0$, $i = 1, ...., g$. For the other one, we easily verify that $\Omega$ is independent of the basis $(\omega_1, .., \omega_g)$. Since the matrix $A$ is invertible, so a base change we can consider $A = I$: $A_{ij} = \delta_{ij}$. Hence $\Omega_{ij} = B_{ij}$, and it still follows from the Proposition1:

$$0 = \int_X \omega_i \wedge \omega_j = \sum_{k=1}^{g} \left( \int_{a_k} \omega_i \int_{b_k} \omega_j - \int_{a_k} \omega_j \int_{b_k} \omega_i \right)$$

$$= \int_{b_i} \omega_j - \int_{b_j} \omega_i$$

Finally, if $v = (v_1, .., v_g) \in \mathbb{R}^g - \{0\}$, then we have:

$$^t v . \operatorname{Im} \Omega . v = \frac{i}{2} \int_X \eta \wedge \overline{\eta} > 0, \text{ when } \eta = \sum_{k=1}^{g} v_k \omega_k$$

∎

## 3 Lattice of periods

Let $X$ be a compact Riemannian surface with two $2g$ fixed simply closed curves which form a basis of the space $H_1(X, \mathbb{Z})$, $(\omega_1, .., \omega_g)$ a basis of the space $\Omega^1(X)$ of holomorphic 1-forms is fixed. The image of the folloing map

$$\begin{aligned} p : H_1(X, \mathbb{Z}) &\longrightarrow \Omega^1(X)^* \\ \gamma &\longrightarrow p(\gamma) \end{aligned}$$

is a lattice $\Lambda$ in $\Omega^1(X)^*$, where $p(\gamma)(\omega) = \int_\gamma \omega$.

**Dfinition 6** *We call $\Lambda$ the lattice of periods. The dual basis $(\omega_1, .., \omega_g)$ identifies the space $\Omega^1(X)^*$ to $\mathbb{C}^g$. As a lattic in the space $\mathbb{C}^g$, $\Lambda = A\mathbb{Z} + B\mathbb{Z}$*

**Remarque 7** *Note that the set $\Lambda$ is a lattice since it comes from the Riemann bilinear relations and the real range of $(A, B)$ is equal $2g$. The Riemann bilinear relations even show that $\Lambda$ is a particular lattice.*

**Dfinition 8** *A divisor on a Riemannian surface is the data of a finite set the points $(P_i, n_i)$, wheited by nonzero inegers. The set of divisors is naturally equipped with a commutative group structure. It is a $\mathbb{Z}$-module generated by $X$. A diviser is called effective if its degree $\sum_i n_i = 0$, and the divisor $D$ is principal if $D = \operatorname{div}(f)$ is given by the pôles and zeros of a meromorphic function $f$.*

**Notation 9** $D = \sum_i n_i P_i$, $\deg D = \sum_i n_i$

## 4 Abel-Jacobi map

Wether $O$ and $P$ are two points of a Riemann compact surface $X$. Two paths $\gamma$ and $\gamma'$ link $O$ to $P$ in $X$ differ only by a factor of $H_1(X, \mathbb{Z})$. In another word: $p(\gamma) = p(\gamma') \pmod{\Lambda}$. For any path $\gamma$ the following map

$$\begin{aligned} u_O : X &\longrightarrow \mathbb{C}^g / \Lambda \\ P &\longrightarrow \left( \int_\gamma \omega_1, ..., \int_\gamma \omega_g \right) \end{aligned}$$

is well defined, but depending on the point $O$. Moreover, for each point $P \in X$ we can associate the divisor $P - O$ of degree zero. A divisior $\text{div}(f)$ associated to a meromophic function $f$ is also of degree zero.

**Dfinition 10** *The set of divisors of degree zero is naturally an Abelian group. We call group of Picard $\text{Pic}_O(X)$ the quotient of divisor group of degree zero by the sub-group of divisors associated to meromorphic functions*

**Proposition 11** *The map $u_O$ extends naturally into a group morphism:*

$$\begin{array}{rcl} u : \text{Pic}_O(X) & \longrightarrow & \mathbb{C}^g/\Lambda \\ \sum_P n_P P & \longrightarrow & \sum_P n_P u_O(P) \end{array}$$

*which does not depend on the point $O$*

**Preuve.** Let's show first the map $u$ is well defined. Wether

$$\text{div}(f) = \sum_P n_P$$

where $f$ is a méromorphic function and we set

$$\omega = \frac{df}{2i\pi f}$$

We note $F_k(z) = \int_O^z \omega_k$ for $k = 1, ...., g$. So Proposition 3 implies

$$\sum \text{Res}\left(F_k \frac{df}{f}\right) = \sum_{j=1}^g \left(\int_{a_j} \omega \int_{b_j} \omega_k - \int_{a_j} \omega_k \int_{b_j} \omega\right)$$

The right side is a linear combination in integers of periods $\int_{a_j} \omega_k$, $\int_{b_j} \omega_k$ as integer, because the periods of the 1-form $\omega$ are integers (Resudue formula). The left side is equal to

$$\sum_P n_P F_k(P)$$

Finally the $k^{th}$ coordonate of the image $u_O(P)$ equals $F_k(P)$. Whether we change the point $O$ in another one $O' \in X$ in another one, then

$$\left(u_O - u'_O\right)\left(\sum_P n_P P\right) = -\sum_P n_P \left(\int_O^{O'} \omega_1, ..., \int_O^{O'} \omega_g\right)$$

But the sum of the right hand is zero, because the degree $\sum_P n_P = 0$ ∎

**Dfinition 12** *The map $u$ defined as above is called the Abel-Jacobi map*

**Thorme 13** *(Abel) The Abel-Jacobi map is injective*

**Preuve.** Whether $D = \sum_P n_P P$ is a divisor of degree zero such that $u(D) = 0$, we will finde a meromorphic function $f$ such that $D = \text{div}(f)$. Indeed we will construct a 1-form

$$\omega = \frac{df}{2i\pi f}$$

Let $\omega$ be a 1-meromorphic form on the surface $S$ with simples pôles in the points $P$ of divisor $D$ with residues $n_P$. Hence once again by Proposition1:

$$\begin{aligned} u(D) &= \sum_P n_P u_O(P) \\ &= \sum \text{Res}(u_O \omega) \\ &= \sum_{j=1}^{g} \left( \int_{a_j} \omega \int_{b_j} \omega_k - \int_{a_j} \omega_k \int_{b_j} \omega \right)_{k=1,\ldots,g} \end{aligned}$$

We will modify $\omega$ so that all its periods will become integers: ∎

**Lemme 14** *Whether $x_1, .., x_g, y_1, .., y_g$ are complexe numbers, then there exists a holomorphic 1-form $\eta$ such that*

$$\int_{a_i} \eta = x_i \quad \text{and} \quad \int_{b_i} \eta = y_i$$

*if and only if*

$$\sum_{k=1}^{g} \left( y_k \int_{a_k} \omega_i - x_k \int_{b_k} \omega_i \right) = 0 \quad i = 1, \ldots, g$$

**Preuve.** As the matrix $A$ is invertible, then the vectors

$$\left( \int_{a_1} \omega_1, \ldots, \int_{a_g} \omega_g \right) \quad i = 1, \ldots, g$$

are linearly independent. Now the following linear map is surjective

$$\begin{aligned} \Phi : \mathbb{C}^{2g} &\longrightarrow \mathbb{C}^g \\ (x_1, .., x_g, y_1, .., y_g) &\longrightarrow \left( \sum_{k=1}^{g} \left( y_k \int_{a_k} \omega_k - x_k \int_{b_k} \omega_i \right) \right)_{i=1,\ldots,g} \end{aligned}$$

So $\dim \ker \Phi = g$. But if $\eta$ is a holomorphic 1-form, $\eta \wedge \omega_i = 0 : i = 1, ..., g$, and then Proposition1 implies

$$\left(\int_{a_1} \eta, ..., \int_{a_g} \eta, \int_{b_1} \eta, ..., \int_{b_g} \eta\right) \in \ker \Phi$$

The lemma follows from that the dimension of the space of
the holomorphic 1-forms is equal to the geneus $g$. Since $u(D) = 0$ in the quotient $\mathbb{C}^g/\Lambda$, then there exists integers $(A_1, .., A_g, B_1, .., B_g)$ such that

$$\sum_{k=1}^{g}\left(\left(\int_{a_k} \omega - B_k\right) \int_{a_k} \omega_i - \left(\int_{a_k} \omega - A_k\right) \int_{a_k} \omega_i\right) \quad i = 1, ..., g$$

So by the lemma above, there exists a holomorphic 1-form $\eta$ such that all the periods of the 1-form $\eta - \omega$ are integers. Hence we can consider that $\omega$ has ineger periods. A primitive of the form between $O$ and $z$ gives the meromorphic function

$$f(z) = \exp\left(2i\pi \int_O^z \omega\right)$$

which is well defined, satisfying $\mathrm{div}(f) = D$ ∎

**Thorme 15** *(Jacobi) The Abel-Jacobi map is injective*

**Preuve.** The map $u$ is a group morphism. So it suffices to show that the image of the map $u$ contains a neighborhood of the point $O$. This will follow from the inverse function theorem: ∎

**Lemme 16** *There exists $g$ distincts points $P_1, .., P_g \in X$ such that any holomorphic 1-form which vanishes in each $P_k$ is identically zero*

**Preuve.** For any point $P \in X$ the sub-space

$$H_P = \{\omega \in \Omega^1(X)^* : \omega(P) = 0\}$$

is of codimension $\leq 1$ in $\Omega^1(X)$. But the intersection

$$\bigcap_{P \in S} H_P$$

is trivial and $\dim \Omega^1(X) = g$. Then there exists points $P_1, .., P_g \in S$ such that

$$H_{P_1} \cap ... \cap H_{P_2} \cap H_{P_g} = 0$$

Let $P_1, .., P_g \in X$ be fixed points as in the lemma with simply connected disjoint local coordonates $(U_i, z_i)$ around these points and $z_i(P_i) = 0$ $i \leq g$. In fact each 1-form $\omega_i$ is written as:

$$\omega_i = \varphi_{ij} dz_j \text{ on } U_j$$

The matrix $(\varphi_{ij})_{1 \leq i,j \leq g}$ is invertible by lemma above.

Consider now the following map

$$F : U_1 \times ... \times U_g \longrightarrow \mathbb{C}^g$$
$$z = (z_1, .., z_g) \longrightarrow (F_1(z), .., F_g(z))$$

such that

$$F_i(z) = \sum_{j=1}^{g} \int_{P_j}^{z_j} \omega_i : i = 1, ..., g$$

The integral

$$\int_{P_j}^{z_j} \omega_i$$

is well defined since each $U_i$ is simply connected. Hence the map $F$ is differentiable in complexe coordonates $z_1, .., z_g$ and the expression of the jacobian matrix is

$$\left(\frac{\partial F_i}{\partial x_j}\right)_{1 \leq i,j \leq g}(P) = (\varphi_{ij}(P))_{1 \leq i,j \leq g}$$

This matrix is invertible in the point $P = (P_1, .., P_g)$. So by the local inverse theorem we have a neighborhood of $F(P) = 0$:

$$W = F(U_1 \times ... \times U_g) \subset \mathbb{C}^g$$

Finally if $\xi \in W$ then there exists points $Q_1, .., Q_g \in \mathbb{C}^g$ such that

$$\left(\sum_{j=1}^{g} \int_{P_j}^{Q_j} \omega_1, ..., \sum_{j=1}^{g} \int_{P_j}^{Q_j} \omega_g\right) = \xi$$

In another wordrs

$$u\left(\sum_{j=1}^{g}(Q_j - P_j)\right) = \xi$$

∎

Summarizing the theorem of Abel-Jacobi:

**Thorme 17** *(Abel-Jacobi)* *The Abel-Jacobi map* $u : Pic(X) \longrightarrow Jac(X) = \mathbb{C}^g/\Lambda$ *is bijective*

Furthermore whether a point $O \in X$ is fixed, we have the following map

$$\begin{aligned} u_O : X &\longrightarrow Jac(X) \\ P &\longrightarrow u(P - O) \end{aligned}$$

When $g = 1$ this map is an isomorphisme. In general it is still:

**Proposition 18** *If the genus* $g \geq 1$, *the map* $u_O : X \longrightarrow Jac(X)$ *is an embedding*

**Preuve.** Since $S$ is compact, it suffices to show that $u_O$ is an injective immersion map. Let's prove firstable $u_O$ is injective. Suppose by contradiction that $u_O(P) = u_O(P')$. So the map $u$ concels on the divisor of degree zero, $P - P'$. This last is the divisor of a meromorphic function $f$. This one has a single pole and a single zero; so it is a map:

$$X \longrightarrow \mathbb{CP}^1$$

of degree one. Thus is absurde since $g \geq 1$. Let's prove that $u_O$ is an immersion map. As in the proof the Abel-Jacobi theorem:

$$d_P u_O(\xi) = (\omega_1(P)(\xi), ..., \omega_g(P)(\xi))$$

The proposition follows again from the local inverse theorem and the next lemma ∎

**Lemme 19** *The holomorphic 1-forms* $(\omega_1, .., \omega_g)$ *have no common zero*

**Preuve.** Once again by contradiction: if a point $P$ is a common zero. According to Riemann-Roch theorem: the dimension of the space of holomorphic functions having more then one simple pole in $P$ equals:

$$\begin{aligned} &\deg u_O - g + 1 + \dim\{\omega \in \Omega^1(X) : \omega(P) = 0\} \\ = \ & 1 - g + 1 + \dim\{\omega \in \Omega^1(X) : \omega(P) = 0\} = 2 \end{aligned}$$

Then there exists a function $f \in X$, which has a unique simple pole in $P$. So it is a map $f : X \longrightarrow \mathbb{CP}^1$ of one degree, when even an absurdity since $g \geq 1$ ∎

**Remarque 20** *Once a point $O \in X$ is fixed we have more generally a map*

$$X^{(g)} = X^g/\mathfrak{S}_g \longrightarrow Jac(X)$$
$$(P_1,..,P_g) \longrightarrow u\left(\sum_{j=1}^{g}(P_j - O)\right)$$

*and $X^{(g)}$ can be provided with an analytical structure. We showed that the map $X^{(g)} \longrightarrow Jac(X)$ is surjective. For reasons of dimensions we can verify that is finite fibers. We can show:*

- *$X^{(g)}$ and $Jac(X)$ are algebraic variety*
- *The map $X^{(g)} \longrightarrow Jac(X)$ is birationnal*